# Optimization model for electric aircraft tow tractors considering operator coalition

Dan-Wen Bao[a], Jia-Yi Zhou[a], Di Kang[b], Zhuo Chen[a]

[a] College of Civil Aviation, Nanjing University of Aeronautics and Astronautics, Nanjing 210016, China

[b] Department of Civil and Environmental Engineering, Rutgers University, Piscataway, NJ, 08854, USA

E-mail address:

Dan-Wen Bao: baodanwen@nuaa.edu.cn

Jia-Yi Zhou: Zhoujiayi_00@nuaa.edu.cn

Di Kang: di.kang@rutgers.edu

Zhuo Chen: zhuochen@nuaa.edu.cn

Corresponding author: Di Kang

E-mail address: di.kang@rutgers.edu

**Abstract:** Horizontal collaboration between operators can save traffic operation costs, a concept that has been particularly validated in the logistics field. Due to the increasing transportation demand and the introduction of charging times for electric vehicles, there is growing pressure on airport ground support services. This study introduces the concept of coalition to airport ground support services. Unlike most studies on horizontal collaboration, the scheduling for airport ground support services needs to particularly consider punctuality issues. In this study, we separately establish electric vehicle scheduling models for operator-separate and operator-cooperated modes and design an algorithm based on the concept of the Adaptive Large Neighborhood Search algorithm, aiming to obtain the electric vehicle scheduling plan with minimum cost and delay time. Furthermore, the study proposes a cost allocation method that considers the degree of sharing among operators to ensure the feasibility of coalition. Finally, we conducted numerical experiments based on actual airport operation data. The experiments verified the effectiveness of the algorithm and the cost allocation method. Compared to solvers, our algorithm can obtain feasible solutions in a shorter time while ensuring that the objective function gap is within 2%. Additionally, the improved cost allocation method is fairer compared to the traditional Shapley method. Numerical experiments also show that coalition can save 15-25% of airport operating costs and 26-39% of delay time, with savings varying based on the sharing parameters. Through quantitative analysis such as sensitivity analysis, the study provides insights into the variation patterns of overall and individual shared utilities and offers suggestions and decision-making mechanisms for the configuration and operation of airport ground operators.

**Keywords:** Horizon collaboration; Airport ground support service; Electric vehicle routing problem; Collaboration; Flight delay; Cost allocation.

## 1 Introduction

As the global civil aviation industry gradually recovers in the post-pandemic era, the demand in the aviation market has experienced significant growth. According to the *Global Outlook for Air Transport* released by the International Air Transport Association (IATA) in 2023, the demand for air travel is expected to double by 2040, with an average annual growth rate of 3.4%. The number of departing passengers is projected to increase from approximately 4 billion in 2019 to 8 billion by the end of the forecast period. However, this growth poses challenges to the efficiency and environmental sustainability of airport ground support services, which require cooperation among various ground support service operators. To reduce carbon emissions, these operators are starting to introduce electric ground service equipment and vehicles to replace traditional fuel-powered fleets. For instance, Beijing Daxing Airport in China has already deployed a ground service fleet composed entirely of electric vehicles.

Airport ground support service operators can be divided into two main categories. The first category consists of the airports themselves, which act as ground support service operators primarily for non-base airline flights. The second category includes base airlines or contractors hired by airlines. The number of these operators depends on the number of service operators introduced by the base airlines at the airport, and they only provide ground support services for the flights of the airlines that hire them. Currently, although multiple operators offer services at airports, their service targets are fixed, and there is no cooperation between them.

To further improve the efficiency of airport ground support services and reduce delays caused by the unavailability of ground service vehicles, this paper proposes a coalition model. Using towing services as an example, the participating operators intend to share tractors. The concept of horizontal cooperation is widely used in urban transportation and logistics distribution and has been proven to offer significant cost savings (Defryn & Sörensen, 2018; Ergun et al., 2007). Distinguishing itself from existing studies, the coalition model for airport ground support services discussed in this paper exhibits the following distinctive traits:

1. Electrification of the fleet: The focus is on a fleet of electric tractors, which consume a large amount of energy during towing and require significant downtime for charging, making scheduling more challenging.

2. Consideration of delay time optimization: Delay time is a critical metric for assessing the service level at airports. This study considers the substantial impact of delay time on operations.

3. Complexity in allocating delay costs: The qualification of delay costs reveals that the impact differs across airports and airlines. For instance, airlines normally bear a greater portion of the costs associated with passenger dissatisfaction. Unlike previous studies where the unit delay costs were often the same for different operators, this paper explores new methods for distributingdelay costs among operators.

To clearly illustrate the problem addressed in this paper, Figure 1 is shown to demonstrate the application of coalition at airports and the potential for resource optimization. Figure 1(a) presents the network of the example, where two service operators, denoted as operator R and operator B, are represented by red and blue markers, respectively. In this example, operator R has 4 flights awaiting service, while operator B has 2 flights, depicted as square nodes. The numerical intervals below the nodes represent the planned start time interval. It is assumed that each operator has only one tractor available for service. Figure 1(b) illustrates a scheduling scheme under the traditional service model (i.e., without coalition), where the trajectory for operator B is depot B-1-2-3-4-depot B, and for operator R is depot R-5-6-depot R. Figure 1(c) presents a scheduling scheme under coalition, where flight 3 originally assigned to operator B is now assigned to operator R. In this scenario, the trajectory for operator B becomes depot B-1-2-4-depot B, while for operator R it becomes depot R-5-3-6-depot R. The intervals above the nodes in Figures (b) and (c) display the calculated service start and service end time intervals, and when the left number exceeds the planned start time interval, it is considered a service delay.

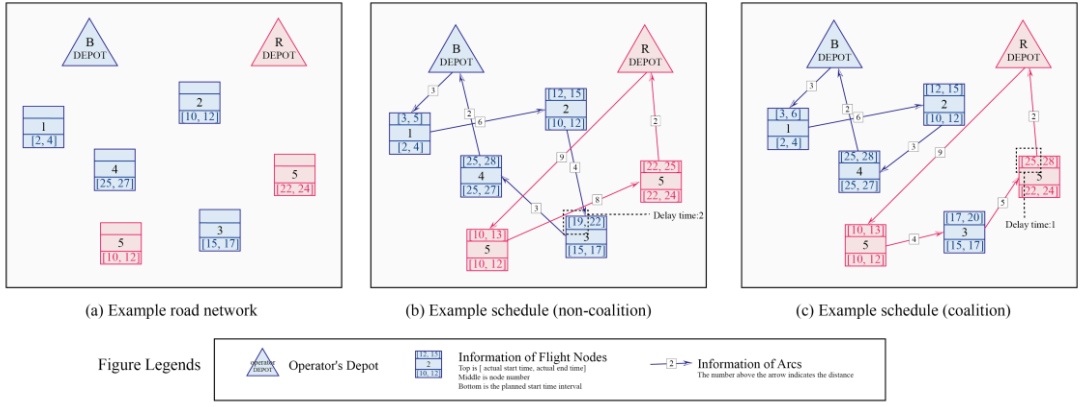

**Fig. 1. Example of problem**

By calculating the distances traveled and arrival times for the two modes described above, the scheduling plan under the segregated operation in this small-scale example involves a total travel distance of 37. Considering the scheduling plan under coalition, the total travel distance is reduced to 34, resulting in a savings of 3.34% in distance traveled. Additionally, under the non-cooperation mode, there is a 2-minute delay at node 3, whereas under cooperation, there is only a 1-minute delay at node 6. Consequently, the total delay time is reduced by 50%. This small-scale example illustrates to some extent the economic and efficiency advantages of coalition, but does not represent an optimal schedule, and this paper aims to maximize the shared utility in real-world scenarios.

This paper addresses the scheduling problem of electric tractors considering a coalition among operators and proposes a cost allocation method based on the degree of resource contribution. The main contributions are:

1. For the first time, a coalition among operators is considered in the context of airport ground operations. Three shared parameters are proposed for scenarios with incomplete sharing: sharing radius, number of shared vehicles, and service priority. A coalition among airport operators is established, focusing on electric tractors, with the aim of optimizing operating costs and delay times.

2. To solve the multi-objective model, a $\xi$ constraint based COADH algorithm is proposed. Additionally, new destroy and repair operators are developed, which are more suitable for solving multi-objective problems considering operator coalition, to quickly and accurately search for high-quality feasible solutions.

3. The paper examines the cost distribution mechanism under conditions of unequal unit delay costs and different sharing intensity among operators. It introduces the measurement of resource contribution and proposes a cost allocation method based on the Shapley value, enhancing the stability of coalition.

4. By conducting a case study, the performance of the algorithm, the effectiveness of the model, and the advantages of the cost allocation method are explored. The impact of different parameters on the model and the effect of individual operator parameters on the stability of the coalition are analyzed, providing data reference and decision-making support for practical applications.

The remainder of the paper is structured as follows: The second section reviews the literature on airport support vehicle scheduling and horizontal cooperation scheduling. The third section introduces two scheduling models for airport electric tractors under operator-separate mode and operator-cooperate mode. The fourth section presents the relevant solution algorithms and the improved cost allocation method proposed in this paper. The fifth section analyzes the impact of various parameters on scheduling results and examines the overall and individual savings achieved through coalition in a case study.

## 2 Literature review

In recent years, **horizon collaboration** has attracted the attention of researchers and practitioners in the transportation field. Horizon collaboration specifically refers to multiple service providers working together at the same level of the supply chain (Cleophas et al., 2019). Previous studies have widely demonstrated the benefits of horizon collaboration across various fields: (Cruijssen et al., 2007) investigated the potential benefits and obstacles of horizon collaboration among logistics companies in Flanders. Their results showed that while service providers strongly believe in the potential benefits of horizon collaboration, the vast majority of companies think that small companies might lose customers or even be pushed out of the market. Aloui et al. conducted a broader systematic literature review on collaborative sustainable transportation (Aloui et al., 2021). Finding mechanisms that make horizon collaboration advantageous for the involved companies is a topic of great interest in the transportation industry (Gansterer & Hartl, 2020).

This section briefly overviews the current state of research on horizon collaboration relevant to our study. For more detailed information, refer to the work by Catherine Cleophas and colleagues (Cleophas et al., 2019), which presents the state of research on various collaboration models, not limited to horizon collaboration. Within the framework of horizon collaboration, hybrid integer programming models are typically developed by extending some of the existing classic vehicle routing models. Algorithms are then designed to solve these models to provide solutions for collaboration decisions(Kang, Li, et al., 2022). For instance, Cruijssen expanded the single-depot vehicle routing problem with time windows to estimate the synergistic benefits of combining outsourcing with horizon collaboration (Cruijssen et al., 2007); Dahl and Derigs developed a dynamic vehicle routing problem model with time windows for pickups, solving it using greedy decoding-based indirect local search (Dahl & Derigs, 2011). By providing this solution to all partners, the decision support system proposed plans for request allocation and routing. Elena Fernández extended the Multi-Depot Vehicle Routing Problem (MDVRP) model, solving the collaborative vehicle routing problem with shared customers using a branch-and-cut algorithm. In this problem, some customers require services from multiple operators, who can decide to transfer part of the demand to their alliance partners (Fernández et al., 2018). E. Angelelli and colleagues considered fairness in their model, adding constraints to control the exchange of workloads between companies to prevent bias towards certain companies (Angelelli et al., 2022). Christof Defryn and colleagues extended the selective vehicle routing problem, allowing alliance partners to set urgency values for their original customers before pooling all requests. When a demand is not served, the urgency value corresponding to the service is added to the objective function to change the importance of the demand (Defryn et al., 2016). Hanan Ouhader considered both economic and ecological impacts, establishing a multi-objective model to study the trade-offs and relationships between multiple objectives under horizon collaboration (Ouhader & El Kyal, 2023). Most of the above studies focus on the logistics field, while vehicle scheduling research in the airport domain more often considers time factors(Kang & Levin, 2021), typically focusing on VRPTW(Norin et al., 2012) or EVRPTW (Bao et al., 2023), but has not yet considered operator cooperation.

   In addition, profit or cost allocation in horizon collaboration is an issue that cannot be overlooked. Research on cost sharing within alliances has proposed several classic game-theory-based methods, such as the Shapley value (Dai & Chen, 2012), providing a fundamental approach to cost allocation in urban transportation collaborations. Some studies have suggested synchronizing cost allocation with routing. For example, Vanovermeire and Sörensen considered the idea of integrating transportation orders in a decentralized manner, proposing that Shapley value-based cost allocation be integrated into the routing problem, rather than first consolidating demands and then allocating costs (Vanovermeire & Sörensen, 2014). Defryn also proposed combining cost allocation with collaborative routing. They modified the Shapley value method based on urgency values to prevent operators from increasing their completion rate by artificially inflating service urgency values (Defryn et al., 2016). Mosleh Amiri and 's research simultaneously considered transferable (cost) and non-transferable (customer coverage) objectives. They proposed a generalized form of the Shapley value with constraints based on the generalized core concept and developed an algorithm to find the optimal cost allocation(Amiri & Farvaresh, 2023). Some scholars have attempted to study fairer allocation methods (Fernández et al., 2016; Gansterer et al., 2021), proposing an optimization model based on the arc routing problem to simulate collaborative operations in truck transportation. The optimization model sets a lower bound for individual operator profits to ensure that all operators benefit from the collaboration.

   In summary, research on ground support vehicle scheduling in the context of horizon collaboration at airports is relatively lacking. While studies on horizon collaboration in urban areas are more abundant, they generally focus on optimizing economic or ecological impacts, making them difficult to directly apply to airports where service time requirements are stringent. Additionally, due to the varying degrees of sharing and differing unit delay costs among operators in this study, existing cost allocation methods are not suitable. Therefore, it is necessary to propose a cost allocation method tailored to the scheduling model in this study.

## 3 Methodology

### 3.1 Problem Description

   The subject of this study is electric tractors that provide towing services for aircraft during taxiing. It is worth noting that the proposed methodology is applicable to other ground service vehicles without loading restrictions. To explore the effects of operational modes at airports on vehicle scheduling outcomes, this paper assumes that all tractors are identical and possess uniform performance characteristics. Within the study airport, multiple operators at the apron are responsible for managing the tractors. These operators typically include various base airlines and airport authorities. Each operator manages a fleet responsible for towing aircraft prior to departure. In the practices of real airports, the service modes utilized by operators can mainly be categorized into two types: operator-separate mode and operator-cooperate mode:

- *Operator-Separate mode.* The base airlines at the airport operate their own fleets of ground service vehicles, dedicated to servicing the flights operated by themself. Base airlines are responsible for covering the operational costs of their fleets and bearing the losses incurred from flight delays caused by ground service. Non-base airlines, on the other hand, rely on the airport to organize ground service. As intermediaries, the airport assumes the operational costs and compensates airlines for delays at a pre-agreed rate.

- *Operator-Cooperate mode.* In this mode, the service range of ground service vehicles (belonging to both the base airlines and the airport) is expanded beyond their own flights, operating based on a pre-established sharing strategy. Participating operators assign a portion of their vehicles to serve a chosen set of flights from other operators as eligible for service and establish specific service priorities (further details in Section 3.4). This model aims to reduce operational expenses and minimize delays for operators.

    To better illustrate the establishment process of the operator-cooperate mode and provide theoretical support for subsequent comparative experiments, we develop two separate models. The first model is for the Operator-Separated Electric Tractors Routing Problem (OS-ETRP), which serves as the basic framework for the traditional operator-separate mode. This model aims to minimize both energy consumption and delay time, representing dual objectives. The second model, designed for the Operator-Cooperated Electric Tractors Routing Problem (OC-ETRP), incorporates the operator-cooperate mode by adding cooperative operation constraints to the first model. Additionally, the second model considers the scenario where vehicle resources are not entirely shared, managed through restrictions on service radius, the number of vehicles being shared, and service priorities. To address the issue of heterogeneous fleets within the same operator, the OC-ETRP introduces the concept of virtual operators. Furthermore, constraints are introduced to regulate the degree of sharing among operators, controlling the extent of shared vehicles. The relationships between these constraints and the degree of sharing are discussed further in subsequent subsections.

    The sets, parameters, and variable descriptions are as follows. Within the time window $[0, t^{LMAX}]$, there is a ground road network $G = (V, A)$ in the airport, where $V$ represents the set of nodes and $A$ denotes the set of arcs. Each individual arc is expressed by $(i, j) \in A, i \in V, j \in V \setminus \{i\}$. The node set $V$ comprises a set of flight nodes $F$, a set of depot nodes $O$, and a set of charging nodes $E$. Each flight node $f \in F$ represents an aircraft parked at a designated aircraft stand, awaiting towing services within a predefined time window as specified in the flight scheduling process. There are multiple operators in consideration, denoted by set $R$, and the set of all vehicles controlled by these operators is denoted by set $K$. For each operator $r \in R$, $F_r \subset F$ represents the set of flight nodes that operator $r$ can serve. The departure depot node for operator $r$ is denoted by $o_r$, the arrival depot node is denoted by $o_r{}'$, and $O = \{o_r \bigcup o_r{}' \mid \forall r \in R\}$.

    For each arc $(i, j) \in A$, $d_{ij}$ represents the distance between node $i$ and node $j$; $t_{ij}$ represents the travel time between node $i$ and node $j$; $h_{ij}$ represents the energy consumption between node $i$ and node $j$.

For each flight node $f \in F$, the parameters are described as follows: $K_f$ represents the set of vehicles available for servicing flight node $f$; $R_f$ represents the set of operators capable of providing service for flight node $f$; $t_f^E$ denotes the earliest start time required for servicing flight node $f$; $t_f^L$ denotes the latest start time; $s_f$ represents the duration of one towing service; and $h_f$ represents the energy consumption of one towing service.

For each vehicle $k \in K$, the parameters are described as follows: the battery capacity is represented by $Q$; Vehicle speed is denoted by $v$; Energy consumption rate per meter is represented by $\tau^e$; Charging rate is denoted by $\tau^c$; and the minimum energy threshold is represented by $\gamma$.

All notations for the OS-ETRP are summarized in Table 1.

**Table 1**

Parameters and variables of model for OS-ETRP.

| | Category | Symbol | Meaning |
|---|---|---|---|
| Set | Graph | $F; f \in F; F_r \subset F$ | Set of flight nodes; Individual of flight node; Flight node set of operator $r$ |
| | | $O; O = \{o_r \bigcup o_r{}' \mid r \in R\}$ | Set of depot nodes; $o_r$, $o_r{}'$ respectively represent the depot node for the operator $r$'s outbound and inbound operations. |
| | | $E; e \in E$ | Set of charging nodes; Individual of charging node |
| | | $V = F \bigcup E \bigcup O$ | Set of all nodes |
| | | $A; (i,j) \in A, i \in V, j \in V \setminus \{i\}$ | Set of arcs; Individual of arc |
| | | $G = (V, A)$ | Graph |
| | Vehicle | $K; k \in K$ | Set of vehicles; Individual of vehicle |
| | | $K_f \subset K$ | Set of vehicles available for flight node $f$ |
| | | $K_r \subset K$ | Fleet of operator $r$ |
| | Operator | $R; r \in R; R_f \subset R$ | Set of operators; Individual of operator; Set of operators available for flight node $f$ |
| Parameter | Flight | $t_f^E$ | Earliest start time for flight $f$ |
| | | $t_f^L$ | Latest start time for flight $f$ |
| | | $t^{LMAX}$ | The maximum value of the time slice |
| | | $s_f$ | Service duration for flight $f$ |

| | | |
|---|---|---|
| | $h_f$ | Service energy consumption for flight $f$ |
| Graph | $d_{ij}, i \in V, j \in V \setminus \{i\}$ | The distance between node $i$ and node $j$ |
| | $t_{ij} = d_{ij} / v, i \in V, j \in V \setminus \{i\}$ | The travel time taken between node $i$ and node $j$ |
| | $h_{ij} = d_{ij} / \tau^e, i \in V, j \in V \setminus \{i\}$ | The energy consumption between node $i$ and node $j$ |
| Vehicle | $v$ | Vehicle speed (assuming consistency) |
| | $\tau^e$ | Power consumption rate (assuming consistency) |
| | $\tau^c$ | Charging rate |
| | $Q$ | Maximum battery capacity |
| | $\gamma$ | Minimum battery level threshold (%) |
| | $n_r$ | The number of vehicles for operator $r$ |
| Cost | $c^e$ | Unit energy cost |
| Variable | $x_{ijk}, (i, j) \in A, k \in K$ | Decision variable, $x_{ijk} = 1$ when vehicle $k$ traverses arc $(i, j)$, and 0 otherwise |
| | $a_{ik}, i \in V, k \in K$ | Decision variable, the time when vehicle $k$ arrives at node $i$ to start servicing/charging |
| | $b_{ik}, i \in V, k \in K$ | Decision variable, the time when vehicle $k$ leaves node $i$ |
| | $y_{ik}, i \in V, k \in K$ | Intermediate variable, remaining energy of vehicle $k$ upon arrival at node $i$ |
| | $t_{fk}^{End}, f \in F, k \in K$ | Intermediate variable, the time when the vehicle $k$ finishes serving the flight node $i$ |
| | $t_{ek}^{End}, e \in E, k \in K$ | Intermediate variable, the time when the vehicle $k$ finishes charging at the charging node $e$ |
| | $s_{ek}, e \in E, k \in K$ | Intermediate variable, the duration of vehicle $k$ charging at the charging node $e$ |
| | $t_f^{DL}, f \in F$ | Intermediate variable, the delay duration of |

The objectives of this paper are as follows: 1) Schedule tasks for electric tractors based on the operator-cooperate mode to complete all aircraft towing services; 2) Ensure sufficient charging time and calculate the timetable for each tractor's arrival and departure; 3) Minimize total operating costs and overall delay time.

## 3.2 Assumptions

This paper focuses on solving the electric tractor scheduling problem under the operator-cooperate mode while identifying the critical parameters influencing total operational costs. To maintain focus within this scope and avoid unnecessary complexity, the following assumptions are adopted:

1. This paper simplifies the ground road network in the airport, focusing solely on three types of nodes: charging nodes, depot nodes, and flight nodes. Road congestion and vehicle conflicts are ignored.

2. This paper assumes uniform performance for all tractors, including travel speed, charging speed, total battery capacity, and vehicle model, while electric tractors are generally non-leveraged and compatible with all flights.

3. Vehicle travel speed remains constant throughout the network. Acceleration, deceleration, and their effects on battery consumption are not accounted for. Therefore, this paper assumes a linear relationship between battery consumption and distance traveled.

4. Electric tractors follow a full charging strategy, meaning they must charge to 100% capacity at charging stations before they leave.

5. This paper assumes that operators establish predetermined service time windows for requests during flight scheduling process. Furthermore, it does not account for changes in estimated departure times resulting from flight delays.

6. The paper does not account for the cascading effects of delays in pre-flight services (such as refueling, catering, and boarding), allowing tractors to start services promptly within the designated time windows.

## 3.3 Operator-Separated Electric Tractors Routing Problem

In the OS-ETRP, operators exclusively serve flights belonging to their affiliated airline company. The model developed for OS-ETRP is based on the classic Electric Vehicle Routing Problem with Time Windows (EVRPTW). It extends the fleet set of the classic Electric Vehicle Routing Problem with Time Windows (EVRPTW) (Raeesi & Zografos, 2020) and defines specific objective functions for airport operations, which are shown in Equations (1) – (2).

$$min\ f_1 = c^e \Sigma_{k \in K} \Sigma_{(i,j) \in A} d_{ij} x_{ijk} \tag{1}$$

$$min\ f_2 = \Sigma_{f \in F} t_f^{DL} \tag{2}$$

The model is a dual-objective optimization model. One objective is to minimize total energy consumption. Since we assume a linear relationship between vehicle energy consumption and travel distance, this objective function is influenced by the distance traveled. In addition, according to the Airport Development Reference Manual (ADRM) published by the International

Air Transport Association (IATA), airport delay is a critical indicator for assessing airport quality. Thus, minimizing service vehicle delay time is established as the second objective function.

The constraints in this optimization model can be categorized into four types: trajectory planning, time window restrictions, energy consumption limitations, and variable constraints.

$$\Sigma_{j\in V\setminus\{o_r\}}x_{o_r jk} = \Sigma_{j\in V\setminus\{o_r\}}x_{jo_r'k} = 1, \forall r\in R, \forall k\in K_r \tag{3}$$

$$\Sigma_{j\in V\setminus\{i\}}x_{jik} - \Sigma_{j\in V\setminus\{i\}}x_{ijk} = 0, \forall i\in F\bigcup E, \forall r\in R, \forall k\in K_r \tag{4}$$

$$\Sigma_{k\in K_r}\Sigma_{j\in V\setminus\{i\}}x_{fjk} = 1, \forall f\in F, \forall r\in R_f \tag{5}$$

Equations (3) – (5) relate to the trajectory planning for tractors. Specifically, equation (3) ensures that each vehicle departs from the depot at the beginning and returns to the depot after servicing all requests. Equation (4) constrains that a tractor must exit a (flight or charging) node once it enters. Equation (5) stipulates that for each flight awaiting service, one and only one tractor controlled by an eligible operator is required to provide service. It is important to note that in OS-ETRP, for all flights $f\in F$ awaiting service, the set of eligible operators $R_f$ contains only one element, which is their original operator.

$$a_{jk} \geq b_{ik} + t_{ij}x_{ijk} - t^{LMAX}(1-x_{ijk}), \forall(i,j)\in A, \forall k\in K \tag{6}$$

$$t_{fk}^{End} = a_{fk} + s_f, \forall f\in F, \forall k\in K_f \tag{7}$$

$$t_f^E \leq a_{fk}, \forall f\in F, \forall k\in K_f \tag{8}$$

$$s_{ek} = (Q - y_{ek})/r_c, \forall e\in E, \forall k\in K \tag{9}$$

$$t_{ek}^{End} = a_{ek} + s_{ek}, \forall e\in E, \forall k\in K \tag{10}$$

$$b_{ik} \geq t_{ik}^{End}, \forall i\in F\bigcup E, \forall k\in K \tag{11}$$

$$t_f^{DL} = \Sigma_{k\in K}(max(0, a_{fk} - t_f^L)), \forall f\in F \tag{12}$$

Equations (6) – (12) represent the time window constraints for tractors. Equation (6) expresses the relationship between the arrival time at the next node in the planned trajectory and the departure time from the previous node. For each tractor, Equation (7) calculates the completion time of service at the flight node. Equation (8) constrains the start time of service for flight node $f$ to be after the required earliest start time. Equation (9) calculates the duration of charging. Equation (10) calculates the completion time for charging. Equation (11) constrains that, for each tractor, its departure time from the flight node must be after the completion time of service or charging. Equation (12) outlines the computation of service delay time, which is a non-linear constraint.

$$y_{jk} \leq y_{fk} - (h_{fj} + h_f)x_{fjk} + Q(1-x_{fjk}), \forall f\in F, j\in V, f\neq j, k\in K_f \tag{13}$$

$$y_{jk} \leq Q - h_{ej}x_{ejk}, \forall e\in E\cup\{o_r\}, j\in V, e\neq j, \forall k\in K \tag{14}$$

$$Q\cdot\gamma \leq y_{ik} \leq Q, \forall i\in V, \forall k\in K \tag{15}$$

Equations (13) – (15) represent the energy consumption limitations for tractors. Equation (13) restricts the energy level upon arrival at the next node after departing from a flight node following service provision. Equation (14) ensures that tractors depart from the depots or charging stations at full battery capacity and calculates the energy level at the next node. Equation (15) restricts the range of vehicle battery levels.

$$x_{ijk}\in\{0,1\}, \quad \forall k\in K, \quad \forall(i,j)\in A \tag{16}$$

$$a_{ik}\in N, \forall i\in V, \forall k\in K \tag{17}$$

$$b_{ik}\in N, \forall i\in V, \forall k\in K \tag{18}$$

$$y_{ik}\in N, \forall i\in V, \forall k\in K \tag{19}$$

Equations (16) to (19) describe the variable constraints of OS-ETRP, where $x_{ijk}$ is a binary variable, and the others are integer variables.

## 3.4 Operator-Cooperated Electric Tractors Routing Problem

In the OC-ETRP, each operator can contribute a portion of their vehicle resources to the coalition and, under some policy restrictions, serve flight node demands from other operators. When an operator allocates a portion of their vehicles to the coalition, two types of fleets exist within the operator: a shared fleet and a non-shared fleet. These two fleets have different sets of serviceable flight nodes. A shared fleet has the capacity to serve a greater number of flight nodes. For example, the shared fleet of operator $A$ can serve flight nodes from other participating operators, e.g., operator $B$, within a certain radius around the $A$'s depot, as well as all flights belonging to the $A$. The non-shared fleet only serve operator $A$'s requests. The shared operation mode transforms the scheduling problem from one of homogeneous fleet to heterogeneous fleet scheduling, significantly increasing the complexity of the model. To address this issue and enhance the applicability of the model developed for OS-ETRP, we propose the concept of virtual operator:

**Virtual Operator** In order to avoid unnecessary complexity, in the operator-cooperate mode, an original operator is virtually divided into a non-shared operator and a shared operator. We only divide operators in the scheduling process, while costs are still calculated as if they belong to the same actual operator. To be more specific, the original operators assign a portion of their tractors to create a shared fleet, forming a new virtual operator participating in the coalition. Other tractors remain within the original operator and are only responsible for the flights assigned to the operator specifically, without being involved in shared flights. To distinguish it from original operators, this newly formed virtual operator is denoted as $r'$ and its shared fleet is represented as $K_{r'}$ (Figure 2).

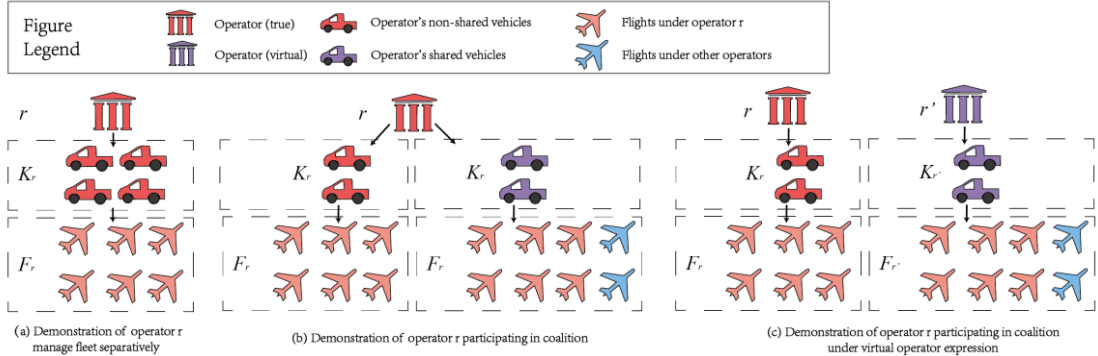

**Fig. 2. Diagram illustrating virtual operator**

In OC-ETRP, prior to scheduling electric tractors, flight nodes need to be assigned among operators to ensure that each flight is serviced by one and only one operator. Before introducing additional constraints, a new decision variable $z_{fr}$ is introduced (equation 20). When it equals 1, it indicates that the flight node $f$ is assigned to the operator $r$; otherwise, it is 0.

$$z_{fr} \in \{0,1\}, \forall f \in F, \forall r \in R_f \qquad (20)$$

Research on horizontal cooperation in logistics has been studied extensively (Fernández et al., 2018), which provides reliable references for this paper. In this section, building upon the model corresponding to OS-ETRP (equations 1-19), while maintaining the objective functions unchanged, we introduce new constraint (equation (21)) and modify equation (5) to equation (22) to establish the model for OC-ETRP.

$$\Sigma_{r \in R_f} z_{fr} = 1, \forall f \in F \tag{21}$$

$$\Sigma_{k \in K_r} \Sigma_{j \in V \setminus \{f\}} x_{fjk} = z_{fr}, \forall f \in F, \forall r \in R_f \tag{22}$$

Equation (21) indicates that for each flight, one and only one operator must be selected to provide service. Equation (22) ensures that the flight nodes are included in the trajectories of the vehicles owned by the operator to which the task is assigned.

For operators with dense distribution of flight in a specific time window, cooperative operations may lead to a shortage of resources due to vehicle being used for external services. To avoid this situation, this paper proposes three key parameters that can influence the degree of sharing among operators in the coalition (Figure 3). Figure 3 illustrates a scenario of cooperative operation between two operators (operator R and operator B). Circular nodes represent depots, triangular nodes represent flight nodes under operator R, and square nodes represent flight nodes under operator B.

Figure 3 (a) illustrates the concept of the service radius. Taking operator R as an example, firstly, it can service all the flight nodes belong to operator R (triangles). Secondly, some of the flight nodes belong to Operator B fall within the service radius of R (hollow squares), so these flight nodes can be serviced by both operator R and operator B. However, the flight nodes belong to operator B outside the service radius of R can only be serviced by operator B (solid squares).

Figure 3 (b) illustrates the concept of the number of vehicles being shared. It shows the possible vehicle trajectories for operator R. Assuming operator R has one shared vehicle and one non-shared vehicle. The trajectory of the non-shared vehicle can only include flight nodes belong to operator R (triangles), while the trajectory of the shared vehicle can include flight nodes that can be serviced by both R and B (triangles or hollow squares).

Figure 3 (c) describes the principle of service priority. Assume that operator R sets a higher priority for its own flight nodes to prioritize servicing its own flights with shared vehicles. When the time window of a flight from operator B overlaps with that of operator R, the vehicle will prioritize servicing Operator R's flights.

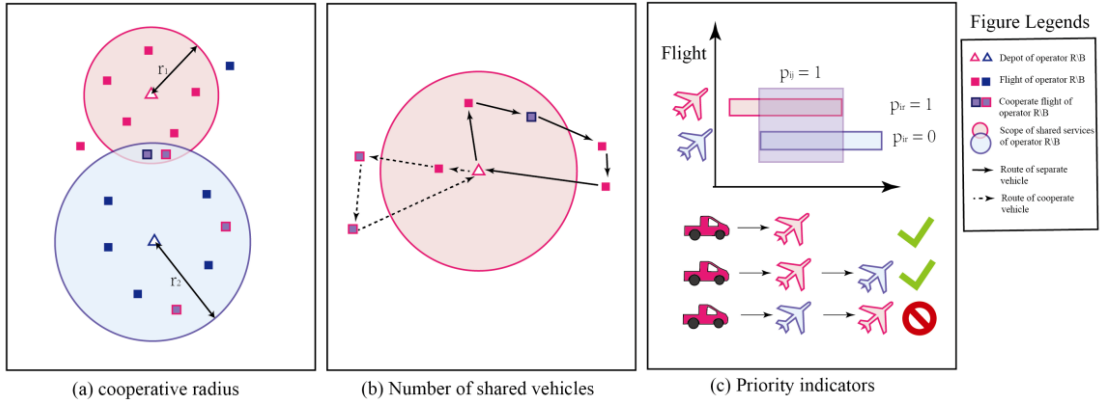

(a) cooperative radius     (b) Number of shared vehicles     (c) Priority indicators

**Fig. 3. Mechanism of the three parameters affecting the degree of sharing.**

Based on the above descriptions, parameters describing the degree of sharing are added in Table 2. These parameters are determined by the operators participating in the shared operations.

**Table 2**

**New parameters for operator-cooperate mode.**

| Catagory | Symbol | Meaning |
|---|---|---|
| Shared parameters | $L_r, r \in R$ | Max service radius of operator $r$ for coalition |

| | |
|---|---|
| $n_{r'}, r \in R$ | Number of shared vehicles provided for coalition of operator $r$ |
| $q_{ij}$ | Flag for whether there is a time window overlap between flight $i$ and flight $j$ |
| $p_{fr}$ | Flight $f$ 's priority for operator $r$ |

Accordingly, three new constraints describing the degree of sharing are added to OC-ETRP:

1) Service radius:

$$F_r = \{f : p_{fr} = 1, f \in F\} \bigcup \{f : d_{fo_r} \leq L_r, f \in F\}, \forall r \in R \qquad (23)$$

Equation (23) restricts that non-native flights serviced by each operator must be within a certain distance from the depot. When $L_r = 0$, it signifies that the operator does not participate in cooperative operations. When the range with $o_r$ as the center and $L_r$ as the radius exceeds the coverage of the airport road network, it means that the operator has no limitation on the service radius, indicating a high degree of sharing.

2) The number of vehicles being shared:

$$|K_{r'}| = n_{r'}, \forall r \in R \qquad (24)$$

Equation (24) restricts the number of vehicles in the shared fleet of an operator. When $n_{r'} = 0$, it indicates that the operator $r$ is not participating in the coalition. Conversely, when $n_r = 0$, it signifies that all vehicles of the operator $r$ are available for sharing, indicating a high level of sharing.

3）Service priority:

$$(z_{ir} - z_{jr})q_{ij}(p_{ir} - p_{jr}) \geq 0, \forall i \in F, j \in F \setminus \{i\}, \forall r \in R \qquad (25)$$

Equation (25) describes a constraint that decides the service priority. When there is a time window overlap between flight $i$ and flight $j$, if an operator $r \in R$ has different priorities for flight $i$ and flight $j$, it is not allowed to assign only lower-priority services to operator $r$. On the contrary, it is feasible to assign all services to operator $r$, assign none, or assign only higher-priority flights to operator $r$. Setting the service priority parameter appropriately can reduce instances of providing services to other operators causing delays for oneself, but it also implies a decrease in the degree of sharing by the operator.

## 3.5 Model Linearization

Because the model includes nonlinear constraints (Equation (12)), linearization is necessary. For ease of representation, we first introduce a new variable $t_{fk}^{DL}$ to denote the delay time of vehicle $k$ at flight node $f$. If the arrival time of vehicle $k$ at the flight node $f$ does not satisfy $a_{fk} - t_f^L \leq 0$, it indicates that the service is not delayed, then:

$$t_{fk}^{DL} = max(0, a_{fk} - t_f^L), \forall f \in F, \forall k \in K_f \qquad (26)$$

Since the objective function is a minimization function, we need only define the lower bound of $t_{fk}^{DL}$. Thus, equation (12) can be replaced by the linearized expressions in equations (27)-(29):

$$t_{fk}^{DL} \geq 0, \forall f \in F, \forall k \in K_f \qquad (27)$$

$$t_{fk}^{DL} \geq a_{fk} - t_f^L, \forall f \in F, \forall k \in K_f \qquad (28)$$

$$t_f^{DL} = \Sigma_{k \in K_f} t_{fk}^{DL}, \forall f \in F \qquad (29)$$

## 4 Solution Algorithm

### 4.1 Handling of multi-objective functions

Chapter 2 introduces two dual-objective optimization models. In this paper, $\varepsilon$-constraint method is adopted to solve the objective functions in order to obtain the Pareto optimal set. The main idea is to select a primary objective function and transform the remaining objective functions into constraints, thereby calculating each sub-optimization objective and obtaining the Pareto optimal set (non-dominated solution set). Since the energy consumption of vehicles varies significantly with the strategy implemented, it is difficult to transform the energy consumption objective into constraints. However, the second objective function (minimizing delay time) has a more straightforward practical significance for constraint conversion. Therefore, we choose to transform the second objective function (equation (2)) into constraint (30):

$$\Sigma_{f \in F} t_f^{DL} < T^{delay} \tag{30}$$

In this equation, $T^{delay}$ represents the maximum delay time of the airport during the operational time period. By adjusting this value, different optimal solutions under various constraints can be obtained, which generate the Pareto optimal solution set. After linearizing and mono-objectifying the models, they become two mixed-integer programming (MILP) models that can be solved by a solver.

### 4.2 Collaborative Operation Adaptive Dispatch Heuristic Algorithm (COADH)

Considering that the models proposed in Section 3 involves numerous variables and constraints, making it difficult to solve using precise algorithm. This section introduces a Collaborative Operation Adaptive Dispatch Heuristic (COADH) algorithm based on the Adaptive Large Neighborhood Search (ALNS) approach. COADH is designed to address the electric tractor scheduling problem under the operator-cooperate mode by modifying certain aspects of the logic within the ALNS framework. ALNS has been widely used in research, particularly for routing and scheduling problems(Pisinger & Ropke, 2019). Ropke and Pisinger first introduced the concept of Variable Neighborhood Search (VNS) (Ropke & Pisinger, 2006), where multiple destroy and repair operators are systematically applied during iterations to improve the current solution. ALNS, on the other hand, can track the success frequency of each destroy and repair operator. Based on this success frequency, it updates the probability of selecting each operator, thus favoring the selection of operators that yield improved performance.

To better solve the model proposed in this paper, COADH modifies traditional ALNS as follows: 1) Since the model revolves around electric tractors, COADH performs separate deletion and addition operations on charging nodes within trajectories. 2) COADH introduces a neighborhood generation method that is more adaptable to dual-objective optimization model and constraints related to coalition (i.e., new destroy and repair operators). Figure 4 depicts the pseudocode of ALNS utilized in this paper, illustrating the general framework of ALNS. The algorithm uses a classic simulated annealing framework as the acceptance criterion for solutions, thereby expanding the search scope and preventing the algorithm from getting stuck in local optima.

---

**Algorithm 1:** Collaborative Operation Adaptive Dispatch Heuristic

**Input:** $T_{max}, \sigma^+, \sigma^-, \sigma, w_1, w_2, w_3, \gamma$
**Output:** x*
initialize feasible solution x;
x* = x;
**while** $T_{max} > 1$ **do**
    *choose destroy number by* $\sigma$
    *choose destroy operator by* $\sigma^-$
    *choose repair operator by* $\sigma^+$
    *destroy x;*
    *repair x;*
    *add charging node, get x';*
    **if** $f_{eval}(x') < f_{eval}(x^*)$ **then**
        x*=x';
        $\sigma^+, \sigma^-, \sigma += w_1$
    **else**
        *generate u in U(0,1);*
        **if** $u \le e^{\frac{f_{eval}(x') - f_{eval}(x^*)}{T_{max}}}$ **then**
            x = x';
            $\sigma^+, \sigma^-, \sigma += w_2$
        **else**
            $\sigma^+, \sigma^-, \sigma += w_3$
        **end**
    **end**
    $T_{max} = T_{max} * \lambda$
**end**

---

**Fig. 4. COADH pseudo-code**

In Figure 4, $T_{max}$ represents the maximum temperature; $\lambda$ represents the rate of temperature reduction; $\sigma^-, \sigma^+, \sigma$ represent the current scores of destroy operator, repair operator, and the number of selections; $w_1, w_2, w_3$ represent the increase in scores under different performances, $x$ represents the current solution, $x'$ represents the new solution, $x^*$ represents the best solution, and $f_{eval}(x)$ represents the fitness function of the solution $x$.

### 4.2.1 Handling for charging station

Previous studies have also utilized ALNS to solve electric vehicle scheduling problems (Bao et al., 2023)(Stenger et al., 2013), In these studies, charging nodes are inserted into the initial solution generation process and are selected for destroy and repair with a certain probability during the iteration of ALNS. This design leads to a fixed number of charging nodes being included in the initial solution generation. In addition, due to the inherent randomness in neighborhood generation, it frequently produces solutions that exceed the battery capacity constraints, which increases the number of iterations and runtime required by the algorithm. In the context of this paper, the setting is an airport where the interruption of airport ground service vehicles due to insufficient energy can result in delays or road congestion.

Therefore, to ensure that the battery constraints are not violated and to give vehicles the opportunity to choose closer charging stations, this paper removes all charging nodes before initiating neighborhood generation. As for the calculation of travel costs, charging nodes are incorporated based on the current arrangement of flights in the trajectories. The steps for vehicle $k$ to insert charging node into trajectory $[o_r, f_1, f_2, ..., f_l, o_r']$ are outlined below:

1) Initialize the battery level $y_{o,k}$ to $Q$;

2) For each subsequent flight node $f \in F$, update $y_{fk}$ based on driving and servicing energy consumption. If $y_{fk} < Q\gamma$, append the charging station $e$ after the flight $f$, and reset the battery level $y_{ek}$ to Q. Repeat this process until we iterate the trajectory to $o_r'$.

3) Update the trajectory timetable and calculate the objective function for this trajectory.

### 4.2.2 Evaluate Function

In ALNS, the evaluate function is a crucial metric for evaluating the quality of solutions. As neighborhood generation allows for infeasible solutions, the evaluate function of the solution in this paper consists of two parts: the first part is the objective function, and the second part is the penalty function resulting from constraint violations. The solutions generated by the algorithm already satisfy most of the constraints. For the remaining potential violations of constraints, we convert them into penalty functions.

The evaluate function for OS-ETRP (denoted as $f_{eval}^{OS}(x)$) considers the penalty function for the first objective function (equation (1)) and the delay constraint (equation (30)), which is mathematically expressed by equation (31). In addition, the evaluate function for OC-ETRP (denoted as $f_{eval}^{OC}(x)$) incorporates penalty functions for the first objective function (equation (1)), the delay constraint (equation (30)), and priority constraint (equation (25)), which is shown in equation (32). In both models, a lower evaluate function value indicates a better solution. Note that when "$f_{eval}^{model}(x)$" does not include a specific superscript, it refers to the evaluate function of the model, i.e., $f_{eval}(x)$ does not specifically denote any particular model's evaluate function.

$$f_{eval}^{OS}(x) = f_1(x) + f_{DL}(x) \tag{31}$$

$$f_{eval}^{OC}(x) = f_1(x) + f_{DL}(x) + f_{PR}(x) \tag{32}$$

In the equations (31) – (32), $f_{DL}(x)$ represents the penalty function for the delay constraint, and $f_{PR}(x)$ represents the penalty function for the priority constraint. These two terms are defined as follows:

$$f_{DL}(x) = c_{eval}^{DL}(\Sigma_{f \in F} t_f^{DL} - T^{delay}) \tag{33}$$

$$f_{PR}(x) = c_{eval}^{PR}(-\Sigma_{r \in R}\Sigma_{i \in F_r}\Sigma_{j \in F_r}(z_{ir} - z_{jr})q_{ij}(p_{ir} - p_{jr})) \tag{34}$$

where $c_{eval}^{DL}$ represents the unit delay cost, and $c_{eval}^{PR}$ represents the unit penalty cost used for evaluating priority constraints. From the perspective of overall scheduling optimization, there is no distinction between unit delay costs among operators. The purpose of these two cost parameters is to standardize the units of measurement between penalty functions and the objective function $f_1$.

### 4.2.3 Destroy Operator

COADH uses the following destroy operators: random removal, worst removal, Shaw removal, distance-based removal, delay-based removal, delay chain removal, and priority removal. The first three destroy operators are well-known in the literature and originate from the work of Ropke and Pisinger (2006)(Ropke & Pisinger, 2006). The fourth to seventh operators are proposed in this paper to adapt to the dual-objective operator-cooperate model. Each operator takes the current solution and parameter $n$ as input, where $n$ defines the number of nodes to be destroyed. The specific descriptions of the destroy operators are as follows:

**Random removal:** randomly removes $n$ nodes.

**Worst removal:** computes the difference between the evaluate function value $f_{eval}$ before and after removing each flight node $f \in F$, sorts them in descending order, and removes the top $D \cdot y^{\chi_{worst}}$ nodes from the list. Here, $D$ represents the size of the list, $y \in [0,1]$ is a uniformly distributed random number, and $\chi_{worst}$ is the parameter for the worst removal operator that measures the impact of changes in objective values on selection. Repeat these steps until $n$ nodes are removed.

**Shaw removal:** Firstly, randomly select a node for removal. Then, in each iteration, randomly select a flight node from the already removed nodes, and sort all remaining nodes in descending order based on their similarity to the selected node $rel(i,j)$. Remove the $\chi_{shaw}$-th node from the list. The similarity between nodes $i$ and $j$ ($i, j \in F$) is calculated as per equation (35). Repeat these steps until $n$ nodes are removed.

$$rel(i, j) = \chi_d \frac{d_{ij}}{\max\limits_{i, j \in F}(d_{ij})} + \chi_e \frac{|e_i - e_j|}{\max\limits_{f \in F}(e_f) - \min\limits_{f \in F}(e_f)} \tag{35}$$

**Travel-cost-based removal:** This operator is a variant of the worst removal operator, where only the first objective function (equation (1)) is considered during the removal process. It calculates the difference between the objective values before ($f_1(x)$) and after ($f_1(x')$) removing each node, sorts them in descending order, and removes the $D \cdot y^{\chi_{worst}}$-th node from the list. Again, $D$ represents the size of the list, $y \in [0,1]$ is a uniformly distributed random number, and $\chi_{worst}$ is the parameter used in the worst removal operator. Repeat these steps until $n$ nodes are removed.

**Delay-based removal:** This operator is a variant of the worst removal operator, which considers the changes in the penalty function $f_{DL}$ related to delay time. It calculates the difference between the penalty values before ($f_{DL}(x)$) and after ($f_{DL}(x')$) removing each node, sorts them in descending order, and removes the $D \cdot y^{\chi_{worst}}$-th node from the list. $D$, $y \in [0,1]$, and $\chi_{worst}$ have the same meaning as defined in the travel-cost-based removal operator. Repeat these steps until $n$ nodes are removed.

**Delay-chain removal:** Due to the transitivity of flight delays, a delay chain is defined as a contiguous segment of the trajectory consisting of two or more consecutively delayed flights, with the length of the delay chain being the number of consecutively delayed flights. The delay chains are sorted in descending order based on their length, and the algorithm selects the longest delay chain from the list. The first flight node in this selected delay chain is then removed. Repeat these steps until $n$ nodes are removed.

**Priority-based removal:** Firstly, randomly select a flight $i$ with priority 0 from the current trajectory. Then, randomly select a flight $j$ with priority 1 from the set of flights served by the operating carrier, ensuring that flight $j$ has overlapping service time with flight $i$. These two flight nodes are then paired together as $T$. Repeat these steps until $n$ nodes are removed. If there are no flights meeting the criteria in the current solution, reselect the operator without changing its score. When selecting and applying this operator, the repair operator is specified as priority-based repair.

### 4.2.4 Repair Operator

The repair operators used in COADH include random insertion, greedy insertion, k-value regret insertion, delay-based insertion, and priority-based insertion. The first three operators are classical operators proposed by Ropke and Pisinger (Ropke & Pisinger, 2006)as well as Molenbruch et al(Molenbruch et al., 2017), while the latter two are new ones proposed for the evaluate function in this paper. Repair operators aim to rearrange the positions of already removed operators in a certain pattern. If a trajectory position meeting the requirements cannot be found, the execution is terminated, and selection is redone based on operator scores. Once all removed flights are reassigned, a new feasible solution is obtained. The specific descriptions of repair operators are as follows:

**Random insertion:** Each flight node $f$ in the set of destroyed nodes is inserted into the solution at any position with equal probability.

**Greedy insertion:** Calculate the increment $f_{eval}$ of inserting each node from the set of destroyed nodes into every position along the trajectory. Select the position with the smallest increment $f_{eval}$ for insertion and repeat this process until all destroyed nodes are reinserted.

***k*th-regret insertion:** Assign a probability proportional to the regret value for each node in the set of destroyed nodes. The regret value is calculated by $r_f^k = \Sigma_{h=2}^{k}(c_f^h - c_f^1)$, where $c_f^1$ is the minimum insertion cost of the flight $f$ in the current solution, $c_f^h$ is its $h$-th ordered insertion cost (in ascending order), and $k$ is the parameter of the operator, typically chosen from 2, 3, or 4. The operator inserts one destroyed node into its optimal position based on probability, updates the regret values of the remaining nodes, and repeats this process until all destroyed nodes are reinserted.

**Delay-based insertion:** Similar to greedy insertion, calculate the increment $f_{DL}$ of inserting each flight node $f$ from the set of destroyed nodes into every position along the trajectory. Select the position with the smallest increment $f_{DL}$ for insertion and repeat this process until all destroyed nodes are reinserted. The purpose of this operator is to minimize delays during execution.

**Priority-based insertion:** When the destroy operator selects removal based on priority, this repair operator is consistently selected. In this scenario, if there are several pairs of flight nodes $T$ in the set of destroyed nodes, for each pair of nodes $T$, swap the positions of the two flight nodes and reinsert them into the trajectory.

### 4.2.5 Adaptive Mechanism

The purpose of the adaptive mechanism is to determine the selection probabilities based on the performance of operators and the number of destroyed nodes in each iteration. For each destroy operator, repair operator, and their corresponding numbers, there are scores $\sigma_i^-, \sigma_i^+, \sigma_i^=$ that continuously change with the number of iterations. The expression for the calculation of the probability $\pi_i^{op}$ of selecting the $i$th operator is shown in equation (36):

$$\pi_i^{op} = \frac{\sigma_i^{op}}{\sum_{m=1}^{M^{op}} \sigma_m^{op}} \tag{36}$$

The subscript "op" in the term $\pi_i^{op}$ represents the destroy operator (-), repair operators (+), or number of nodes to destroy (=). $M^{op}$ represents the total number of operators available for random selection in the adaptive mechanism. In particular, since the priority-based insertion operator is not included in the random selection process, $M^+ = 4$.

In each iteration, COADH updates the score of operators based on the solutions generated by the operators by adding $w$. The value of $w$ is determined by equation (37):

$$w = \begin{cases} w_1, & \text{better solution and accepted} \\ w_2, & \text{not better solution but accepted} \\ w_3, & \text{do not accepcted} \end{cases} \tag{37}$$

## 4.3 Cost Allocation Among Operators

The cost allocation method affects the stability of **coalition**. Shapley value method is considered the general allocation method for horizontal logistics cooperation (Ahari et al., 2024). The Shapley value considers the value contribution of operators to all possible (sub)coalitions, thus being entirely based on the cooperative productivity of operators. The basic formula for the Shapley value is shown in equation (38):

$$\phi_r = \sum_{S \subseteq R, \, \{r\}} \frac{|S|! \cdot (|R| - |S| - 1)!}{|R|!} \cdot (c(S \cup \{r\}) - c(S)) \tag{38}$$

where $\phi_r$ represents the cost allocated to operator $r$; $R$ is the set of all operators; $S$ denotes any subset of operators that does not include operator $r$; $|S|$ and $|R|$ denote the number of operators in the sets $S$ and $R$, respectively; $c(S)$ represents the total cost of the set of operators $S$.

The key to calculating the allocated costs is to compute the total cost $c(S)$ for any subset S. To calculate the total cost of the coalition to be allocated, it is essential to clarify the cost structure of coalition, which consists of travelling cost and delay cost:

**Travelling cost:** since each tractor has a travelling cost, the allocable travelling cost for operator $r$ is the sum of the travelling costs for each trajectory belonging to the corresponding virtual operator $r$.

**Delay cost:** unlike travelling costs, delay costs are borne by the actual operators and are independent of the allocated trajectories. The calculation involves summing the delays for all flights under an operator. It is challenging to separate the delay costs between virtual operators included in coalition and actual operators not included in coalition. Therefore, we propose a delay cost allocation method based on resource contributions $\rho_r$.

The calculation of the allocated cost $c(S)$ for the set of operators $S$ is shown in equation (39):

$$c(S) = \Sigma_{r \in S} (c^E \Sigma_{k \in K_r} \Sigma_{(i,j) \in A} d_{ij} x_{ijk} + \rho_r c_r^{DL} \Sigma_{i \in F_r} t_f^{DL}) \tag{39}$$

where $c^E$ represents the travelling cost per unit distance; $c_r^{DL}$ represents the delay cost per unit for operator $r$, and $\rho_r$ represents the resource contribution of operator $r$. Resource contribution depends on three parameters that influence the degree of sharing as discussed in Section 3.4, including service radius, the number of vehicles being shared, and service priority.

The service radius $L_r, r \in R$ directly affects the number of serviceable flight nodes. In this paper, we assess the contribution of each operator $r$ based on the size of the serviceable flight nodes managed by the shared fleet of operator $r$, the number of shared tractors, and the priority strategy, represented as $\rho_r^1$, $\rho_r^2$, and $\rho_r^3$. The formulas are shown in equations (40-42).

$$\rho_r^1 = \frac{\mid F_{r'} \mid / \mid F_r \mid}{\underset{s \in R}{max}(\mid F_{s'} \mid / \mid F_s \mid)} \tag{40}$$

$$\rho_r^2 = \frac{\mid K_{r'} \mid /(\mid K_r \mid + \mid K_{r'} \mid)}{\underset{s \in R}{max}(\mid K_{s'} \mid /(\mid K_s \mid + \mid K_{s'} \mid))} \tag{41}$$

$$\rho_r^3 = \sqrt{\frac{\sum_{f \in F_r}(p_{fr} - \sum_{s \in F_r} p_{fr} / \mid F_r \mid)^2}{\mid F_r \mid}} \tag{42}$$

In equation (36), $\mid F_{r'} \mid$ represents the number of serviceable flight nodes managed by the shared fleet of operator $r$, and $\mid F_r \mid$ represents the total number of flights managed by operator $r$. In equation (37), $\mid K_{r'} \mid$ represents the number of shared tractors provided by operator $r$, and $\mid K_r \mid$ represents the total number of tractors operated by operator $r$ excluding shared tractors. To eliminate discrepancies in numerical magnitudes between scores, we integrate and normalize individual scores to obtain the resource contribution of operator $r$, denoted as $\rho_r$, with the calculation shown in equations (43-44).

$$\rho'_r = \rho_r^1 + \rho_r^2 + \rho_r^3 \tag{43}$$

$$\rho_r = \frac{\rho'_r - \underset{s \in R}{min}(\rho'_s)}{\underset{s \in R}{max}(\rho'_s) - \underset{s \in R}{min}(\rho'_s)} \tag{44}$$

The preceding description covers the calculation of coalition costs and the allocation of coalition costs among operators. However, for an individual operator, costs comprise both the shared coalition cost $\phi_r$ among operators and those incurred due to independent operation. The independent operating costs include the travelling cost of non-shared tractors and the delay cost not accounted for in coalition costs. Therefore, the final calculation for the costs borne by an operator is given by equation (45).

$$c(r) = \phi_r + c_e \Sigma_{k \in K_r} \Sigma_{(i,j) \in A} d_{ij} x_{ijk} + (1 - \rho_r) c_r^{DL} \Sigma_{i \in F_r} t_f^{DL} \qquad (45)$$

# 5 Numerical Experiments

This section evaluates the performance of the proposed models and the effectiveness of algorithms, as well as the impact of various model parameters. Additionally, we also discuss the cost allocation method and coalition strategies from the perspective of individual operators' interests. The experiments in this paper were conducted on a 64-bit Windows computer with 12th Gen Intel(R) Core(TM) i5-12400F 2.50 GHz processor and 16.0 GB RAM, and the MILP model was solved using Gurobi 12.6.0 software to obtain precise solutions.

The effectiveness of the algorithm is evaluated by comparing it with the features of the precise Pareto optimal solution set. The comparison includes the average values of the two objective functions in the solution set and the computational time. The experiments in this paper were conducted on a 64-bit Windows computer with 12th Gen Intel(R) Core(TM) i5-12400F 2.50 GHz processor and 16.0 GB RAM. The Pareto optimal solution set used for comparison was obtained by solving the MILP model using Gurobi 12.6.0 software.

The effectiveness of the model is evaluated by comparing the operator-cooperated model with the operator-separated model. As described in Section 3, the operator-cooperated model includes equations (1) - (4), (6) - (11), and (13) - (30), while the operator-separated model includes equations (1) - (11), (13) - (19), and (26) - (30).

The remaining of this section are as follows: Section 5.1 introduces the background of the case study. Scenarios are created based on the real road network and flight schedule of Guangzhou Baiyun Airport. In Section 5.2, we describe the tuning process of the algorithm and show the performance of different operators in the algorithm. In Section 5.3, we compare the solution generated by COADH with the exact solution solved by Gurobi to study the accuracy and computational efficiency of COADH. In Section 5.4, we investigate the resources saved by the operator-separated model compared to the operator-cooperated model, propose the concept of shared utility, and test the impact of various important parameters on the results. In Section 5.5, we propose the concept of individual shared utility for the stability of coalition, evaluate the cost allocation method among individuals, and assess the impact of individual shared parameters on individual shared utility.

## 5.1 Background and Scenario Setup

This paper selects the outbound flights from Guangzhou Baiyun Airport (CAN) on February 17, 2024, as the case study. As the largest airport in China, CAN handled a passenger throughput of 631.86 million in 2023, indicating its representativeness among large and busy airports. The full-day flight schedule comprises 758 flight nodes, with a total segment length of 19802.73 meters in the ground road network (Figure 5). The time range of the full-day flight schedule spans from 0 to 1,440 minutes, covering a total of 758 flights operated by 59 airlines. These airlines are categorized into three operators: operator 1 represents airline CZ, operator 2 handles towing services for airlines ZH, HU, MU, JD and GJ, and operator 3 provides services for other non-base airlines using airport service fleet, including all the remaining airlines like FD, ET, AK, etc.. Further details about the operators are provided in Table 3.

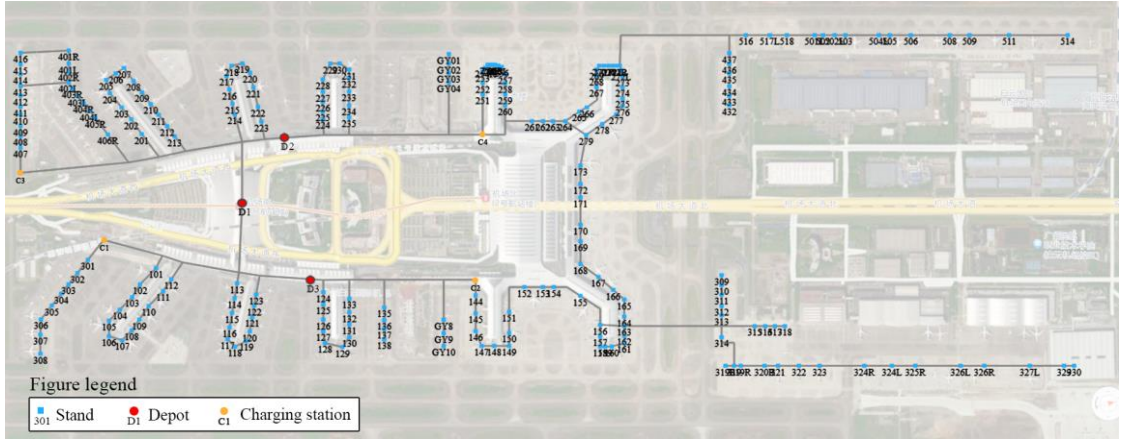

**Fig. 5. Example road network map**

**Table 3**

Table of operators' information.

|  | operator 1 | operator2 | operator3 |
|---|---|---|---|
| Number of flights | 382 | 182 | 194 |
| Number of tractors | 4 | 4 | 4 |
| Max service radius | 1.5km | 1.5km | 1.5km |
| number of shared vehicles | 4 | 4 | 4 |
| Variance of service priority | 0.1 | 0.1 | 0.1 |
| Unit cost of delays (CNY/min/flight) | 8 | 8 | 4 |

The parameters of the tractors are based on data obtained from real airport tractor manufacturers, as detailed in Table 4. Given that tractors from various manufacturers exhibit differing performance characteristics, a sensitivity analysis for the tractor parameters is conducted in Section 5.4.

**Table 4**

**Information of electric tractor.**

| parameter | value |
|---|---|
| Battery capacity (kwh) | 50 |
| Electricity consumption rate (kwh/km) | 0.5 |
| Charging rate (kwh/min) | 0.8 |
| Minimum power threshold (%) | 0.3 |
| Service duration (min) | 3 |
| Mobility rate (km/h) | 10 |

To evaluate the algorithm's performance across varying flight scales on a fixed road network, three scenarios are devised based on flight demand, ranging from minimal to maximal inclusion. The flight settings for each scenario are as follows:

- Scenario A: Encompasses peak-hour flights, totaling 56 flights, involving three operators.
- Scenario B: Encompasses flights during peak hours, as well as the hour preceding and succeeding, totaling 148 flights, involving three operators.
- Scenario C: Encompasses full-day flights, totaling 758 flights, with participation from three operators.

These scenarios serve to assess the algorithm's accuracy in calculations and computation time, alongside evaluating the cost allocation method's effectiveness across different flight demands. Additionally, Scenario B is used for parameter sensitivity analysis. This choice is motivated by the fact that peak-hour scheduling reflects overall efficiency, and extending the timeframe by an hour before and after mitigates bias arising from tractors being in service or charging states during peak hours, thus enhancing the representativeness of Scenario B.

## 5.2 Analysis of ALNS Algorithm

This section presents results related to the parameter tuning process and the effectiveness of the proposed destroy and repair operators. In order to determine the optimal combination of parameters, we conducted tests on parameters $T_{max}$ and $\gamma$. Specifically, we tested $T_{max}$ within the set $\{10^3, 10^4, 10^5, 10^6\}$ and $\gamma$ within the set $\{0.001, 0.005, 0.01, 0.05, 0.1, 0.5\}$. The results show that both of these parameters have an impact on the quality of solutions and computation time. As $T_{max}$ increases, the algorithm performs better. When $T_{max} = 10^4$, the evaluate function reduces by 16% compared to using $T_{max} = 10^3$. However, further increases in $T_{max}$ do not lead to further improvements. Compared to taking the largest value 0.5, when $\gamma = 0.1$, a reduction of up to 30% in the objective function is noticeable. As $\gamma$ continues to decrease, the reduction in the evaluate function becomes smaller and smaller. When $\gamma$ drops below 0.01, no further improvements are observed. Therefore, in the experiments of this paper, we set $T_{max} = 10^4, \gamma = 0.01$. Due to similar data sizes and operating environments, other parameters are based on previous research results(Bao et al., 2023; Levin & Kang, 2023). The summary of all algorithm parameters is presented in Table 5.

**Table 5**

Parameter values for the COADH algorithm .

| Parameter | Value |
| --- | --- |
| $T_{max}$ | $10^4$ |
| $\gamma$ | 0.01 |
| $(w_1, w_2, w_3)$ | (30,18,12) |
| $\sigma_i^{op}\ (initial)$ | 50 |
| $\chi_{worst}$ | 6 |
| $\chi_{shaw}$ | 3 |
| $n$ | [2,4,6] |
| $c_{eval}^{DL}$ | 1 |
| $c_{eval}^{PR}$ | 50 |

Figure 6 illustrates the average score of the destroy operator applied in operator-separated and operator-cooperated models across three different scenarios A, B, and C. In the COADH algorithm, operators with superior performance contribute to the algorithm consistently obtaining improved solutions, consequently yielding higher score values. By comparing the performance rankings of operators, it can be observed that in scenario A, the worst , travel-cost-based, and delay-based operators perform relatively well in both models, consistently ranking within the top three, albeit with slight variations across different scenarios. Conversely, the delay chain, and priority-based operators exhibit worse performance, with significant gaps compared to the top three scores. As the road network and the flight demand increase in scenarios, the performance of worst and Shaw operators improves quickly, while the performance of the remaining operators remains relatively stable. It is worth noting that random operator's performance is not as bad as expected.

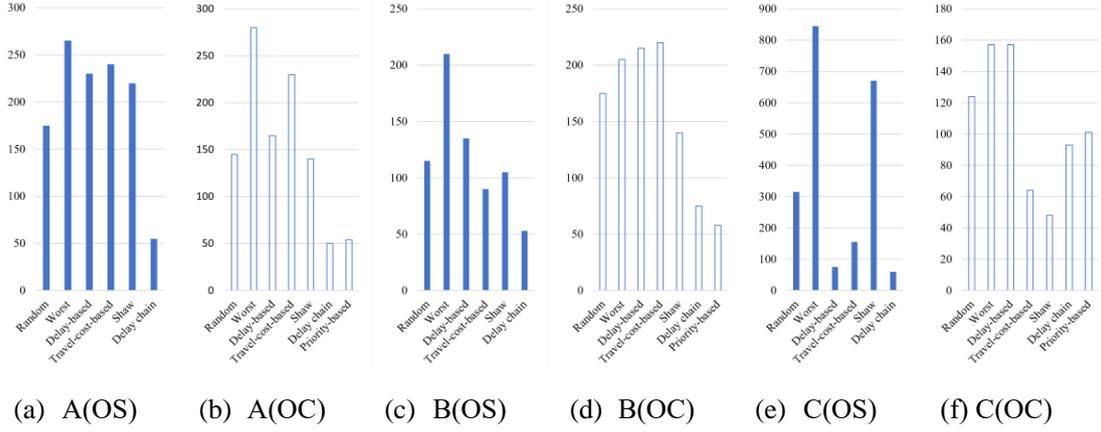

**Fig. 6. Average scores of the destroy operator**

Figure 7 illustrates the average scores of repair operators applied to scenarios of different scales (A, B, C) in both operator-separated and operator-cooperated models. Following the same principle, the greedy and 2th-regret operator performs the best, followed by delay-base and 2th-regret. Conversely, random and priority-based operators perform poorly. Although priority-based repair operator is not a random selection operator, its score aligns with priority-based removal operator. Notably, the 2th-regret operator performs best in smaller scales (e.g., scenario A), but its score gradually decreases with increasing scale. Conversely, the delay-base operator performs better in larger scenarios, especially in OS model.

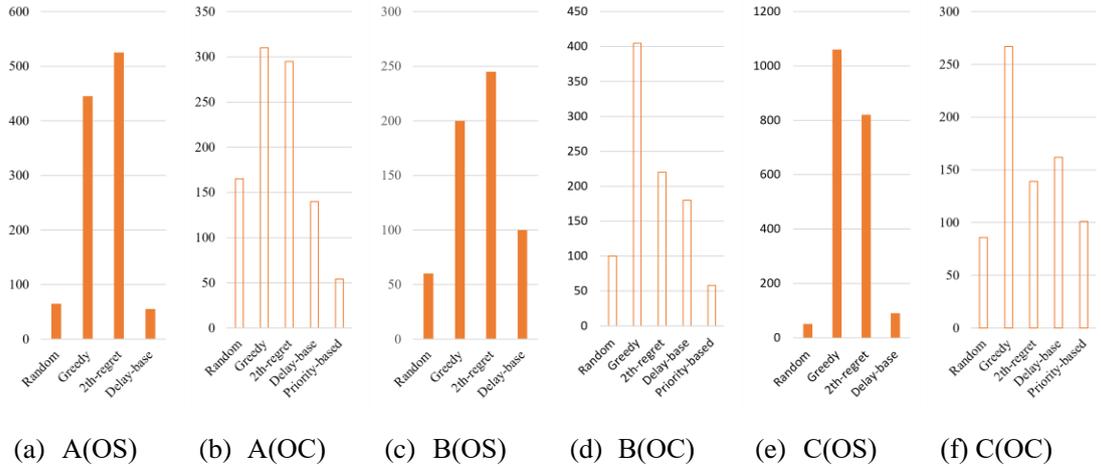

**Fig. 7. Average scores of the repair operator**

## 5.3 Analysis of ALNS Effectiveness

In this section, we present the computation time and accuracy of the ALNS algorithm as applied to scenarios A, B, and C, respectively. Figure 8 shows the variation of the evaluate function with the number of iterations. As depicted in Figure 8, overall, the ALNS algorithm converges, with the convergence speed decreasing as the flight demand increases. Scenario A reaches convergence in about 100 iterations, while scenario C requires around 800 iterations to converge.

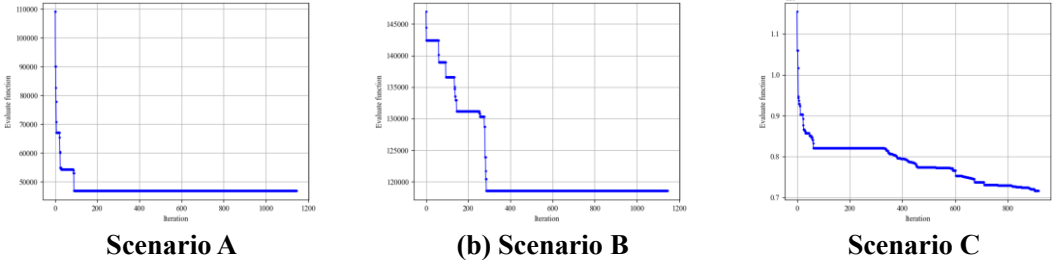

| Scenario A | (b) Scenario B | Scenario C |

**Fig. 8. COADH evaluate function with respect to the number of iterations.**

A more detailed comparison of computational accuracy and computational time is presented in Table 6. The difference in computational accuracy for the two objective functions, as calculated by the average difference in Pareto solutions obtained by COADH and Gurobi, is expressed as a percentage (objective gap column). A higher value in this column indicates higher algorithmic accuracy and time efficiency. Overall, the discrepancy between COADH and Gurobi objectives is below about 2% , which meets the accuracy requirements for practical use.

In comparing the operator-separated model to the operator-cooperated model, COADH demonstrates superior accuracy in optimizing operating costs for the operator-separated model. Conversely, for minimizing delay time, COADH shows greater accuracy in optimizing the operator-cooperated model. In terms of time efficiency, COADH shows significant advantages. In cases where Gurobi can solve the precise solutions, COADH can conserve roughly 90% of computational time, with this savings escalating as the flight demand increases. Additionally, even when Gurobi fails to produce a solution within the fixed time window, COADH can still solve results quickly.

**Table 6**

Comparison of the results solved by COADH and Gurobi.

| TYPE | Model | Scenario A | | | Scenario B | | | Scenario C | | |
|---|---|---|---|---|---|---|---|---|---|---|
| | | COAD H | GUROB I | objectiv e gap | COAD H | GUROB I | objectiv e gap | COAD H | GUROB I | objectiv e gap |
| COST | OS | 49361 | 48941 | -0.9% | 181584 | 179283 | -1.3% | 770903 | * | * |
| | OC | 46834 | 46324 | -1.1% | 140135 | 138974 | -0.8% | 713396 | * | * |
| DELA Y | OS | 23 | 24 | 4.2% | 73 | 72 | -1.4% | 465 | * | * |
| | OC | 31 | 30 | -3.3% | 46 | 45 | -2.2% | 364 | * | * |
| TIME | OS | 17 | 258 | 93.4% | 68 | 721 | 90.6% | 2598 | * | * |
| | OC | 29 | 267 | 89.1% | 105 | 736 | 85.7% | 3240 | * | * |

## 5.4 Experimental Results

### 5.4.1 Shared Utility

To measure the impact of coalition on airport operations, this section introduces the concept of shared utility. Equation (46) describes the calculation of shared utility $\Delta f(L_r, n_{r'}, p_{fr})_i$ :

$$\Delta f(\overline{L_r}, \overline{n_{r'}}, \overline{p_{fr}})_i = \frac{\overline{f_i(S_{os})} - \overline{f_i(S_{oc})}}{\overline{f_i(S_{os})}}, \forall i \in \{1,2\}, f \in F, r \in R \qquad (46)$$

In this equation, $\Delta f(L_r, n_{r'}, p_{fr})_i$ represents the shared utility under the parameters $(\overline{L_r}, \overline{n_{r'}}, \overline{p_{fr}})$ , which denotes the percentage savings of the $i$th objective function of airport towing services with parameters $(\overline{L_r}, \overline{n_{r'}}, \overline{p_{fr}})$ . $\overline{f_i(S_{os})}$ represents the average value of the $i$th objective function of the Pareto solution set obtained by solving the OS-ETRP under certain scenario (i.e., scenarios A, B, or C), while $\overline{f_i(S_{oc})}$ represents the average value of the $i$th objective function of the Pareto solution set obtained by solving the OC-ETRP. According to the description of model objective functions in Section 3, shared utility is categorized into two types based on the type of savings: **shared utility for travelling cost and shared utility for delay time**. In addition, shared utility is categorized into **overall shared utility** and **individual shared utility** based on the subject of savings. The definition of individual shared utility can be found in Section 5.5.

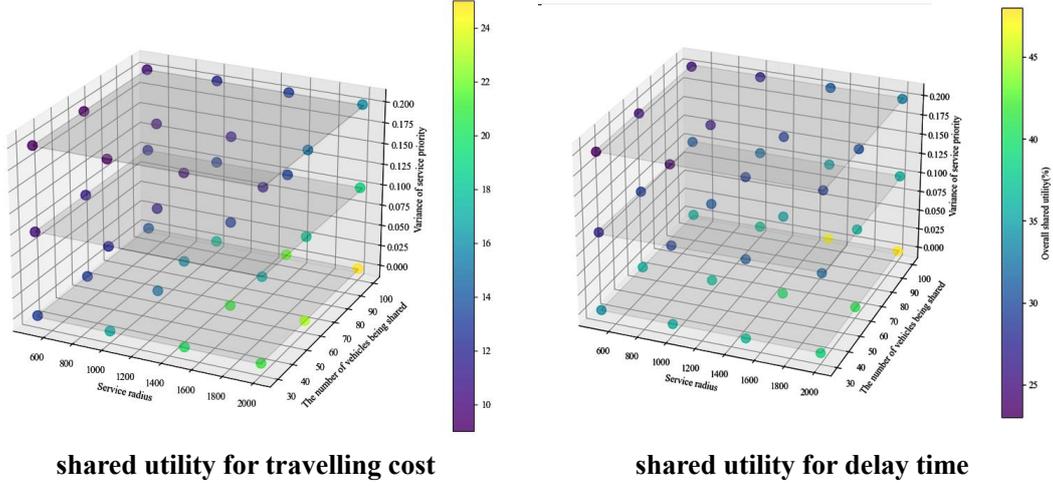

shared utility for travelling cost               shared utility for delay time

**Fig. 9. Impact of coalition on objective functions**

Figure 9 illustrates the overall shared utility of towing services affected by coalition. In general, compared to operator-separated model, operator-cooperated model can save approximately **15-25%** of tractor operating costs and reduce delay time by around **26-39%.** The shared utility is influenced by factors such as the service radius of different operators $L_r$, the number of vehicles being shared $n_{r'}$, and the standard deviation of service priorities $p_{fr}$. With a constant average priority standard deviation, larger service radius and higher proportions of shared vehicles lead to greater overall shared utility. Furthermore, the standard deviation of priorities reflects the extent to which operators differentiate between their own flights and those of other companies. Figure 9 also indicates that a larger standard deviation of priorities results in smaller overall shared utility.

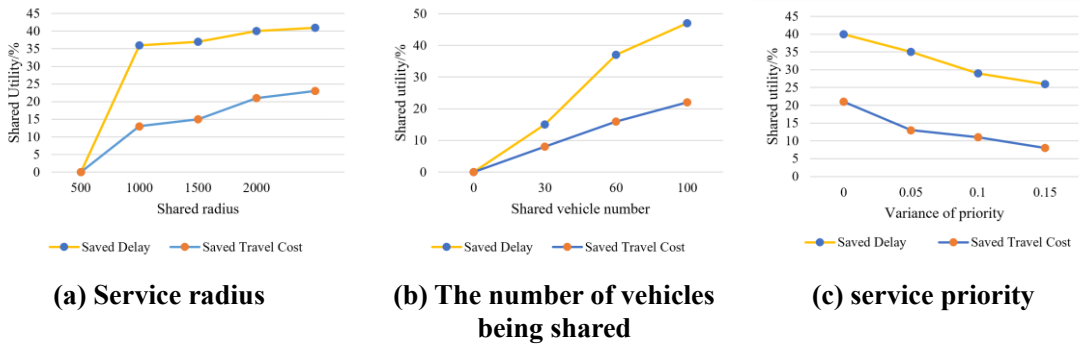

        **(a) Service radius**            **(b) The number of vehicles**        **(c) service priority**
                               **being shared**

**Fig. 10. The impact of shared parameters on shared utility**

Figure 10 provides a more detailed depiction of the influence of three parameters on the objective functions. We established the benchmark for this experiment, where $(\overline{L_r}, \overline{n_{r'}}, \overline{p_{fr}})$ take the values of (1500, 100%, 0). When testing the sensitivity of shared parameter effects on shared utility, we only adjust the range of one parameter at a time, ranging from minimum to maximum influence, while keeping the other parameters constant. In this experiment, it is assumed that all operators adopt comparable levels of shared parameters, with the goal of investigating the influence of overall sharing intensity on operating costs and delay times.

Figures 10 (a) and 10 (b) illustrate the effects of service radius and the number of vehicles being shared on overall shared utility. Both the shared utility for operating cost and shared utility for delay time increase with the increase in service radius and the number of vehicles being shared. Moreover, the shared utility for delay time is larger than the shared utility for operating cost. When the sharing intensity reaches its maximum, the shared utility for operating cost is approximately 20%, while the shared utility for delay time is approximately 40%. Regarding the number of vehicles being shared, the utility increment remains relatively stable, but the increment of the service radius depends on the distribution of the depots and flight nodes, showing a distinct turning point.

Figure 10(c) demonstrates the variation of objective functions as the variance of service priority for each operator changes from 0 to 0.2. As the variance increases, the shared utilities for delay time and operating cost decrease from 40% and 23% to 28% and 11%, respectively. This is because an increase in the service priority variance indicates that operators tend to prioritize their own flights more, thereby reducing the sharing intensity.

From the perspective of a horizontal comparative analysis of different parameters, adjusting the contribution of a shared parameter to the degree of sharing from the lowest to the highest results in an increase in economic sharing utility of 41%, 47%, and 9% for the service radius, number of vehicles, and service priority, respectively. The corresponding increases in efficiency sharing utility are 23%, 22%, and 5%. This indicates that the service priority among the three shared parameters has the smallest impact on shared utility, while the service radius contributes more to the shared utility for operating cost, and the number of vehicles being shared contributes more to the shared utility for delay time.

**5.4.2 Sensitivity Analysis**

Given that certain parameters may vary across different scenarios and are not constants (for example, the power consumption rate depends on the tractors used by the operator), this section conducts analysis on other important parameters that may affect the results of the model. Example B is selected for experimentation, with the shared parameters of the OC-ETRP set to (1500, 100%, 0). This section explores the average values of the objective functions of the Pareto optimal solution sets under different parameters. The discrepancy between the outcomes of the operator-separated model and the operator-cooperated model serves to indicate the utility savings resulting from coalition.

Electric vehicles using charging strategies have different types of charging rates, typically categorized as slow charging and fast charging (Kang, Hu, et al., 2022). The time required for a vehicle to charge from 0 to full battery varies from 30 minutes to 120 minutes. By adjusting the maximum charging time (calculated as battery capacity divided by charging rate, in benchmark mentioned earlier, it is 60 min), within a range of $[30,120]$ minutes, with intervals of 30.

Table 7 presents the average objective function values of the Pareto optimal solution set. Regarding operating costs, the impact of charging time on vehicle operating costs is not significant for operator-separated model, whereas for the operator-cooperated model, operating costs increase by 9.7% as the charging time increases from 30 min to 120 min. Meanwhile, the shared utility for operating cost decreases from 21.8% to 16.4%. In terms of delay time, increasing the charging time leads to a significant increase in delay time for both operator-separated model and operator-cooperated model, with a more significant effect observed in the operator-separated mode. When the charging time increases from 30 min to 120 min, the delay time in operator-separated model increases by 622%, whereas in operator-cooperated model, it increases by only 537%. This difference may be attributed to the limited number of available vehicles per flight in the operator-separated model, consequently reducing scheduling flexibility. Therefore, the shared utility for delay time increases with the increase in charging time. When the charging time reaches 120 minutes, coalition can save 42% of delay time. The aforementioned impacts of charging time on the results can be attributed to faster charging speeds allowing vehicles to have the opportunity to select flights closer within service radius, rather than being more constrained by time limitations. Additionally, the increased availability of idle vehicles at the same time also increases, ensuring more flights are serviced, thereby reducing delay time.

**Table 7**
Impact of charging time on the objective function.

| Full charge time/ min | Dis-os | Dis-oc | $\Delta f_1$ | Del-os | Del-oc | $\Delta f_2$ |
|---|---|---|---|---|---|---|
| 30 | 137536.1 | 107517.0 | 21.8% | 16.7 | 10.8 | 35.1% |
| 60 | 140135.0 | 108880.8 | 22.3% | 46.0 | 28.7 | 37.6% |
| 90 | 138605.7 | 111338.7 | 19.7% | 81.7 | 49.5 | 39.4% |
| 120 | 141216.8 | 117963.1 | 16.5% | 120.7 | 69.1 | 42.7% |

Because electric vehicles of different types have varying battery ranges, it is important to investigate how these ranges affect operations. Range can be measured as the ratio of battery capacity to the vehicle's energy consumption rate, holding the energy consumption rate constant. Based on the existing types of electric vehicles on the market, the battery capacity is set to fluctuate between 60 kWh and 100 kWh, with intervals of 20 kWh. To ensure that the charging time remains unchanged, the charging rate also varies accordingly.

Table 8 presents the average objective function values of the Pareto optimal solution sets for two models. From the perspective of operating costs, whether for operator-separated model or operator-cooperated model, as the vehicle's range improves, there is a slight reduction in the vehicle's operating costs. Higher range capabilities can decrease the number of times vehicles need to go to charging stations, thus reducing a portion of detour distances. The shared utility for operating cost remains nearly unchanged, staying at around 22%. Regarding delay time, an enhanced vehicle range significantly reduces delay time, with a more significant decrease in operator-separated model compared to operator-cooperated model. This is because the reduced charging frequency results in fewer vehicles being in a charging state at any given time, which implies an increase in available vehicles. The $\Delta f_i$ column reveals that as vehicle range increases, the economic sharing utility by coalition essentially remain constant, while the shared utility for delay time decreases.

**Table 8**
Impact of vehicle range on the objective function.

| Battery capacity / kwh | Dis-os | Dis-oc | $\Delta f_1$ | Del-os | Del-oc | $\Delta f_2$ |
|---|---|---|---|---|---|---|

| | | | | | |
|---|---|---|---|---|---|
| 60 | 142749.2 | 110589.4 | 22.5% | 73.0 | 44.8 | 38.6% |
| 80 | 140135.0 | 108880.8 | 22.3% | 46.0 | 28.7 | 37.6% |
| 100 | 138152.1 | 108837.9 | 21.2% | 19.6 | 13.2 | 32.8% |

At different airports, operators have varying numbers of vehicles. The average number of flights served by each vehicle is calculated by dividing the total number of flights by the total number of vehicles. By adjusting the number of vehicles owned by operators, we study the number of vehicles that should be introduced. This parameter is adjusted within the range of [10,16], with intervals of 2 (in benchmark mentioned earlier, it is 12).

Table 9 displays the average objective function values and shared utility (in the $\Delta f_i$ column) for the Pareto optimal solution sets of the two models. Concerning operating costs, as the average number of flights served by vehicles increases, both operator-separated model and operator-cooperated model experience slight reductions in operating costs. Additionally, operator-cooperated model exhibits a higher rate of decrease compared to operator-separated model. This is primarily because scheduling flights on the same trajectory can save the distance vehicles travel from and to the depot, and operator-cooperated model, with a broader selection range of vehicles, can arrange shorter trajectories.

Regarding delay time, both models show a rapid increase in delay time with the increase in the average number of flights served by each vehicle, reaching a maximum overall utility saving of up to 41.9%. According to Table 9, with other parameters unchanged, it can be ensured that the average number of flights served by each vehicle over a 3-hour period remains around 12-14 (1 hour before and after peak hours). Once the average number of flights served by each vehicle surpasses 14, the shared utility for delay time does not increase rapidly. Instead, this situation often results in extended delay times, which runs counter to the airport's service objectives.

**Table 9**

Impact of average flights served per vehicle on the objective function.

| Average number of flights serviced per vehicle | Dis-os | Dis-oc | $\Delta f_1$ | Del-os | Del-oc | $\Delta f_2$ |
|---|---|---|---|---|---|---|
| 10 | 142841.1 | 112860.4 | 21.0% | 2.9 | 2.6 | 10.0% |
| 12 | 140135.0 | 108880.8 | 22.3% | 46.0 | 28.7 | 37.6% |
| 14 | 137952.9 | 106410.6 | 22.9% | 58.3 | 35.7 | 38.7% |
| 16 | 137087.1 | 105060.6 | 23.4% | 99.9 | 61.5 | 38.4% |

The location where vehicles are parked can have a certain impact on the results. We describe the distribution of depots by calculating the distance from the centroid of all flight nodes to the centroid of the depots. By using k-means clustering depots' coordinates, we obtain the coordinates of the centroid and adjust the distribution of different depots to change the centroid position. Table 10 displays the average objective function values and shared utility for the Pareto optimal solution sets of two models. It can be observed that the depot location has little influence on the two objective functions of different models, especially on delay time. However, in terms of shared utility, the closer the depot is to the edge, the higher the cost savings.

**Table 10**

Impact of depot location on the objective function.

| Distance of the garage from the center of gravity of the flight node/ m | Dis-os | Dis-oc | $\Delta f_1$ | Del-os | Del-oc | $\Delta f_2$ |
|---|---|---|---|---|---|---|

| | | | | | |
|---|---|---|---|---|---|
| 0 | 141154.0 | 107523.1 | 23.8% | 45.1 | 28.1 | 37.7% |
| 200 | 140135.0 | 108880.8 | 22.3% | 46.0 | 28.7 | 37.6% |
| 400 | 138772.7 | 112193.8 | 19.2% | 48.9 | 35.5 | 27.5% |
| 600 | 147295.7 | 109670.0 | 25.5% | 43.7 | 30.2 | 30.9% |

## 5.5 Individual Shared Utility

Section 5.4 validated the existence of overall sharing utility. However, in the operator-cooperated mode, whether individual operators have shared utility has a decisive impact on the stability of coalition. Therefore, this section analyzes the shared utility of individual operators. Similar to the overall sharing utility mentioned in Section 5.4, this section defines the individual shared utility of operators, as calculated in Equation (47).

$$\phi_r^{SAVE} = \phi_r^{OS} - \phi_r^{OC}, \forall r \in R \tag{47}$$

where $\phi_r^{OS}$ represents the cost for operator $r$ in operator-separated model, and $\phi_r^{OC}$

represents the cost allocated to operator $r$ in operator-cooperated mode. When $\phi_r^{SAVE}$ is positive, it indicates that joining coalition has saved costs for operator $r$, with larger values indicating better cost savings. Conversely, when $\phi_r^{SAVE}$ is negative, it suggests that joining coalition has increased costs for operator $r$, deteriorating the stability of coalition. In such cases, the airport may operate less efficiently. This section aims to identify cost allocation methods that maximize individual shared utility while ensuring the stability of coalition as much as possible. It also explores the impact of parameters on individual shared utility, providing references for corresponding service providers in the strategy formulation process.

To verify the cost allocation method proposed in this paper (Section 4.3), the cost allocation for the three scenarios mentioned in Section 5.1 is carried out and compared with the traditional Shapley value allocation method. Each solution set contains four solutions, each solution containing three operators. We compute a total of 12 individual shared utilities under each senario, and statistically analyze the number of operators with individual shared utility less than 0 among all Pareto optimal solution sets. results is shown in Table 11. Since the improved Shapley method considers the differences in the degree of sharing between operators, it is more suitable for the coalition proposed in this article. In all three scenarios, new allocation method reduces the ratio of more than 50% of individuals with negative shared utility a result of 0 to 1 individual with negative shared utility.

**Table 11**

The number of operators with negative individual shared utility after cost allocation.

| Senario | Senario A | Senario B | Senario C |
|---|---|---|---|
| Traditional shapley | 6 | 8 | 9 |
| Improved shapley | 0 | 1 | 1 |

By examining the cost allocation results for scenarios B and C in Figure 11, the operator with negative individual shared utility is operator 3 in the solution sets with longer delay time. This operator is characterized by a lower unit delay cost and tight vehicle resources. The reason for the violation of cost allocation constraints is that the coalition, in order to solve the resource scarcity issue of operator 3 itself, exchanged for higher unit delay time from other operators. Although the total delay time is reduced, the total delay cost increases, and according to the allocation rules, more delay costs are allocated to operator 3, thus resulting in a violation of constraints.

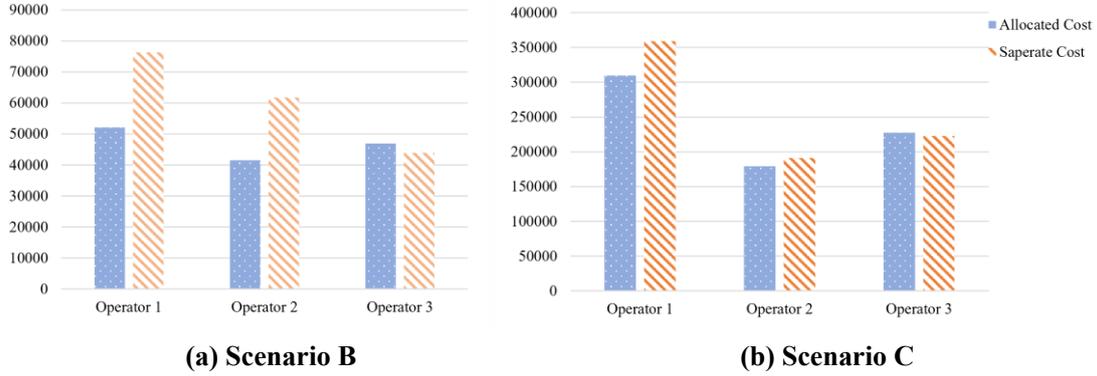

**(a) Scenario B**   **(b) Scenario C**

**Fig.11. Allocation result with negative shared utility**

However, the unit delay cost depends on the characteristics of the operators. For airline operators, the delay cost depends on factors such as the size of the airline, the customer base, and the characteristics of the flights. For airport operators, the delay cost depends on the compensation ratio agreed with the airlines, the size of the airport, and other attributes. These attributes, like vehicles, are difficult to change in a short period of time. Therefore, we attempt to observe changes by individually adjusting the sharing parameters for each operator. Figure 12 shows the changes in cost allocation resulting from adjusting the sharing parameters of operator 3. The main vertical axis in the bar chart represents shared utility $\phi_r^{SAVE}$, while the secondary vertical axis in the line chart represents delay time.

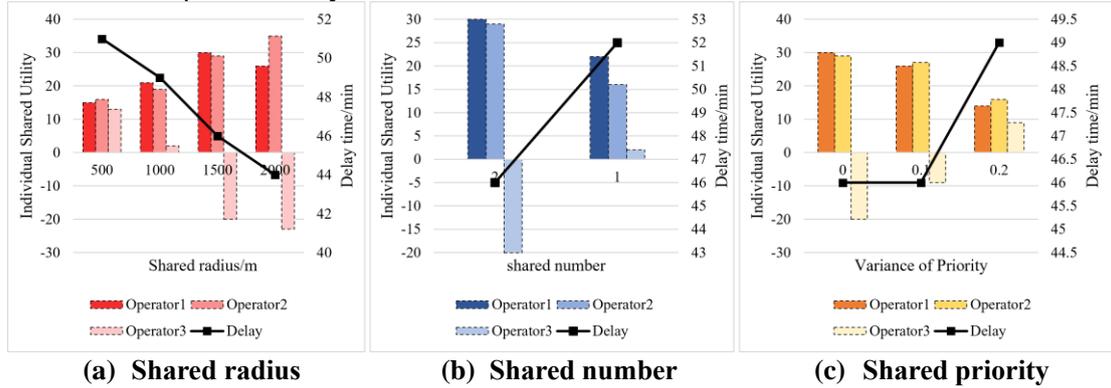

**(a) Shared radius**   **(b) Shared number**   **(c) Shared priority**

**Fig. 12. Impact of Individual Sharing Parameters on Cost Allocation**

Figure 12 shows the changes in individual shared utility with the variation of operator 3's sharing parameters. When the lengths of the bars are similar, it indicates a fairer cost distribution. As seen in Figure 12, adjusting the sharing parameters to reduce the sharing intensity of operator 3—specifically by decreasing the service radius, reducing the number of vehicles being shared, and increasing its service priority—can improve the fairness of cost allocation to some extent and eliminate instances of negative shared utility.

Specifically, when the service radius of operator 3 is reduced by 1,000 km, the number of shared vehicles is decreased to 1, and the priority variance is increased to 0.2, the issue of unfair cost allocation is completely resolved. However, as indicated by the line chart, reducing the sharing intensity slightly increases the total cost, with the number of vehicles being shared having a significant impact on the cost.

When operator 3 takes on more of its own flight tasks, it reduces the resources that other operators need to allocate to operator 3, thereby decreasing the likelihood of increased delay times for operators with high unit delay cost (operator 1 and operator 2) and reducing operator 3's cost allocation in coalition. This prevents cost losses due to the allocation method.

## 5.6 discussion

The numerical experiments provide several management insights. Firstly, the analysis of shared utility reveals that for airports with multiple ground service operators, adopting coalition can save travel costs and improve service efficiency. Operators can set reasonable sharing parameters based on the airport configuration. The data indicates that shared utility increases with the sharing intensity and is most significantly influenced by the number of vehicles being shared, while the impact of priority settings is minimal. This provides a reference for operators when adjusting parameters.

Secondly, the sensitivity analysis offers insights into scheduling results and shared utility under different parameters. Higher shared utility is observed in scenarios with longer vehicle charging times, lower endurance, and more flights served by each vehicle. This suggests that airports with slow charging, short vehicle range, or tight vehicle resources should be encouraged to implement coalition. In the test of scenario B, the coalition cost is lowest when each vehicle serves 12 flights, which can serve as a reference point.

Lastly, the results from the perspective of individual shared utility confirm the effectiveness of the proposed cost allocation method. However, operators with low unit delay time (typically airports) may still experience negative individual shared utility. To address this problem, reducing the sharing intensity of operators with low unit costs can prevent losses for some operators after joining the coalition, though it will increase the overall operating cost. In practice, operators can assess their unit delay costs and adjust parameters according to their specific circumstances and the insights presented here. Alternatively, they can utilize the model and analytical approach outlined in this paper to determine the optimal parameter settings.

## 6 Conclusion

This paper focuses on electric tractors used in airport ground support services. It builds a multi-objective model for electric vehicle routing problem with time windows that incorporates shared operations. By adding sharing variables and constraints to the EVRPTW model, we establish a vehicle scheduling model for operator coalition aimed at minimizing travel cost and delay time. Since the model considers delay time as an objective, it is formulated as a nonlinear programming model. In the context of operator coalition, we innovatively propose three sharing parameters to describe their impact on the sharing intensity and their constraint formulas in the mathematical model. Subsequently, we linearize the delay time calculation equation and use a $\varepsilon$ - constraint method to handle the delay objective function. Then, based on the ALNS concept, we design the COADH algorithm, developing destroy and repair operators more suited for the OC-ETRP. Finally, through case study, we validate the effectiveness of the algorithm and the savings from coalition. We also test the impact of sharing parameters and vehicle parameters on the objective function and shared utility.

From a practical perspective, this research not only provides a vehicle scheduling method with the coalition but also develops a pre-decision mechanism for airport ground service operators to assess the shared utility before implementing coalition. By considering the impact of sharing parameters and certain vehicle parameters on overall and individual sharing utility, the research can also help operators choose the optimal coalition strategy and facility configuration. Additionally, it can provide numerical references for extra compensation contracts among operators to facilitate the formation of coalition agreements, reduce airport operating costs, and improve airport operational efficiency.

Future research could further optimize the model by considering more realistic vehicle energy consumption models or studying vehicle types with loading restrictions. Additional optimization objectives, such as minimizing emissions or scheduling fairness, could also be considered. Finally, to address the issue of individual shared utility, the model could include constraints that ensure this utility is greater than zero, establishing a model that simultaneously considers vehicle scheduling and cost allocation.